\documentclass{article}
\usepackage{latexsym,amsmath,amssymb}

\begin{document}
	
\title{Volume Product of convex  bodies}
\author{R.Anantharaman\\ 
Retired from SUNY College , Old Westbury  NY 11568-0210\\	
rajan.anantharaman@gmail.com} 
\date{ }

\maketitle

\textbf{To  the memory  of Som Naimpally,  teacher and friend}

\begin{abstract}
	
	In this    expository paper we   discuss  the  volume product  P(K) of convex bodies K  in $R^n$; this is the product of volumes of K  and  its polar  K*. The Blaschke- Santalo inequalities  state that always  $ P(K) \le  P(B_2)$  and $ P(B_1)\le P(K)$ .  Here the Closed unit Ball in  norm of $l_2$  is denoted by $B_2$  and like wise for the $l_1$ unit ball. The upper bound is classical, due   to Santalo in general and Blaschke and  Mahler much  earlier in 1930 s for dimensions 2 and 3. The lower bound is open for general  K. However   the result  of  Gordon, Meyer  and Riesner  says that  the  class of  zonoids  K  attain  the lower bound.  There is Bianchi and Kelly' s proof of the upper  bound in general,involving Paley-Wiener Theorem  , as generalized to $R^n$  by Stein .   For the lower bound, there is the result of  of Kim and Zvavitch on its  stability under perturbations of  unconditional  K .  Further there are  results on " functional" versions of Blashke -Santalo Inequality . In this context we discuss  the  newer concept of   "Polar  f* " for certain type of functions  f and its relevance here and mention  as a "functional" example  a result by S.Artstein, B.Klartag and V.Milman .  We mention a later one  by Huang  and Ai- Jun Li  and discuss  Ball's Inequality  for unconditional bodies,another  strengthening of the Blaschke-Santalo  inequality .
\end{abstract}

\underline{MR Classifications:} Primary 52A20,  52A40, Secondary 42A05 

\section{Introduction  and Preliminaries}\label{S:intro}	
	
   \textbf{1.1}  We discuss  the volume product P(K) of convex bodies in $R^n$ . This is the product of volumes of K  and its polar, $K^*$, also  called Mahler  Product  in [19]. The Blaschke-Santalo inequality  and its reverse, have two parts combined here for briefness ; the  upper bound is proved ( see below)  but lower bound is open in this generality :-- 
   
   \begin{equation} \label{E:cong2}
      P(B_1) \le  P(K) \le  P( B_2) 
      \end{equation}
 Here $B_2$ and $B_1$   are the closed unit balls  in the norms  of  $l_2$  and $l_1$.  We  note the equality   $P(K)  =  P(K ^*) $  thanks   to the Bi Polar Theorem.   Further  we only  need  in Section 3 the actual value of the lower bound  namely $\frac{4^n}{ n!} $ ( see  [4]).  This last is seen from the volume of the cube = $ 2^n$ and that of $B_1$  from the fact there are $2^n$ " Octants"  and in each the volume is $ \frac{1}{ n!}$.  We  note  in passing from  [3] [4] the curious fact that the volume of the Euclidean ball in $R^n$ ,   (also) becomes rather small for large n  . In words of [4] :--" This is the first indication  that our low  dimensional intuition may lead us astray".  See also [3]  about the need to  correct intuition(  when needed,  as above) and continue until the next ( shock).    In  this  expository  we do  not deal directly with this  high dimensional  "pathology" .    \\
    
    Returning to Volume Product   Inequalities,    in words of [19 ] 'These     show   how roundest   or pointiest    a generak convex body  K can be "    .The upper  bound was known to  Blaschke and Mahler  in the 1930s  or earlier  for dimensions 2 and 3  see [6]. The general  case was proved by Santalo   in 1949; that  the upper bound is assumed   (only) at $B_2$  or an ellipsoid was proved by Petty in 1982.

    Let us npw  describe the results .\\
    
   \underline{ Sec 2 }   due to Bianchi  and Kelly[6] deals with proof of The Blashcke-Santalo Inequality and the ingenious  use of Paley-Wiener  Theorem  [11 [15]] as generalized by Stein[17] .The authors tool is a functional called $\rho(K)$   .Another such one  namely  $\eta(K)$   is  also used   and verifies/reproved   conjecture by Siegel[]  for 2  cases of K.\\

      \underline{Sec3)} This   due to Gordon, Meyer and Reissner  [9] concerns the lower bound    achieved in the class of zonoids K . The methods are classical  Convexity  and in this sense self contained except for a formula for volume of Zonoids [16]\\

   \underline{Sec 4)} due to Kim and Zvavitch[12]    consider perturbations of   unconditional bodies; they  obtain a "Stability Version" of the lower bound in Eq(1), Blashske Santalo Inequality  . Namely   suppose we know that the result holds for a certain class of convex bodies(in this case Unconditional ones). Let us perturb this    a little.  Then we may hope the result to hold in this new  class.  We  discuss one such result .The tool  to "perturb" is the Banach Mazur  " distance " now adapted to Convex bodies. \\
   
  \underline{ Sec 5 } by Huang and Li  [10] on the other hand  deals with the Functional Form of   and other related classical inequalities.    
   Following 10]  among the several we only discuss an  inequality  that they call   Ball's .   We see in [4] the earlier  result namely the Prekopa-Leindler Inequality that is functional form of Blasche=Santalo Inequality. In [10]  [1] [2]  we see the tool   Polar of a non negative (measurable ) function  ie   generalizing Polar of  Convex sets. Since this concept looks at first  non intuitive we   give a chart, spend some time on it , give  examples and " confirmation"  from[1] [2] and hope it is of interest .
   
    We conclude   by   discussing two inequalities,  (31)   and (32) below, the first due to [1]  [2] implies   Blashchke   Santalos Inequality; the second   from[10 ] Ball's \\

    \underline{	Examples}  These illustrate  results esp  in Sec3  and 5.  Return to Ineq(1) of both  Blaschke-Santalo and its reverse:---

   In the Exs we  deal with  inequalities from Sec 3  and  5 and   look for following:--
   
   (A) Verify Inequality  ( and equality if   it is true   )
   
   (B )Identify Polar set  ( recalled  immediately)
   
   (C)   Verify the  Blashchke   Santalo  Inequality (1)

   Let us note that the    integral in  " quadratic"  inequality  (30)   Sec 5  is harder to compute  and we only consider easier exs
   \\

   \textbf{note:} The \underline{Polar} \textbf{} $A ^0 $ of A  is defined   below  under \underline{Basic notions}   .Let us remark here that the verification of the Polar set it is  enough to verify the required inequality  for extreme points in set A. This   method  is useful in all the examples:\\

    \underline{Tools} For needed Harmonic  Analysis we refer to [11] [16] [18  [21].
   For Convexity and Convex geometry  [4] [8] [13  [14] and[17].  In particular the last 3 contain many more relevant results we donot pursue.  
   
   \underline{Disclaimer}  We  do not treat  many important related  results here. For instance the  earlier but deep  result of Bourgain and Milman on the lower bound ( see references for ex[5],[6],[10][12],[19]) using theory of co type and many  others following or improving- this  one. Let us also merely refer to    [13]  [14] for  the important   results on  the Busemann- Petty  Problem and solution[ 13] [14].
   
   We refer to [ 19 ] [20]for  insights   on   this topic; [17] is a veritable  dictionary on it Vol.Product is in ch10  .  Let us refer   again  to the   Bourgain memorial paper [5] by   Ball  for  related  results  many more  than what we  can consider here.\
   
   Finally let us recall \underline{basic   notions} (for ex[7])  and use these without quoting them leaving   other definitions at the  appropriate places :--
   Throughout we consider Convex bodies K in  $  R^n$ ; also symmetric( most the time)
   The \underline {polar} K* of   K is the set
   
   \[
 (  y in R^n   :   |< x,y>| \le 1  for every x  in K   )
    \]
     
    ie. it is the   "dual ball"  of the  dual space of the Banach space whose unit ball is K.

     We have then  \underline{ the dual norm} 
     
     \begin{equation}\label {E:conj}
      ||y||_{K^*}  =sup  ((x, y): x\epsilon  K ) for   y\epsilon  R^n
      \end{equation}
     
      and called classically  [7] as\underline{ the support function }, $h_K (. )$ of K  .

       Also the \underline{ Minkowski gauge function }  $|| . ||_K$ of K is defined  by
    
    \begin{equation}\label{E:conj}   
     ||x||_K  = inf (  t>0 :  x \epsilon t K)
     \end{equation}  
     
      we  have
   $ h_K(y) = ||y||_K* $    We   use   both: the  other  " support function " definition  is easier to handle \\

 \section{ Upper bound  [6] . Tool  Paley-Wiener-Stein Theorems }

In this Section we  discuss  results from [6]  proving the   upper bound in Blaschke Santalo Inequality ,Eq(1)  in general,  for  convex , balanced bodies  K. 

Let us record some facts, or  tools   from Harmonic Analysis ( [ 11  [16]  [18][21]])that we need  below  for Theorem  2.2  : 

(a) First the \underline{Fourier Transform  or FT}   of  an F in $L^2$ (or $L^1$)   is   defined by 

 $ \hat{F}(\xi ) $   =  $\int F(x)  e^{-2\pi  x.\xi} dx$  for$\xi\epsilon R^n$   and integral over $R^n$
 
 (b)The \underline{Inverse  FT} is $ F( x )    =\int \hat{F}(\xi)  e^{2\pi  x.\xi} dx$ for$ x\epsilon R^n$   and integral over $R^n$.
 For  general conditions  on validity of $ L^1 $ Inversion we refer to above references.Let us   stress  that in our  cases  next   the function F is in $L^2$    or has support K. 
 
(c) The next tool is\underline{ Plancherel's Theorem:} The FT as an operator on $ L^2$  is an   onto isometry.

(d) As mentioned before  the  extension by Stein to $R^n$  of the classical Paley-Wiener Theorem  is the major tool  used. Let us state it now:--\\

 \underline{Paley Wiener  Theorem and its extension  by Stein}

The  important notion here is  "  Function is of Exponential   type "  

Let us recall  the    Theorem in 1  dimensions 

\underline{Theorem [11 ] [156]}Suppose f is  
in $L^2 (R)$ .Then f is the FT of a  function vanishing off the interval I=$[-M ,  M ] $  iff f(x) is  the restriction to R of an entire function  f(z) of \underline{exponential type M  }.ie  it satisfies the following estimate:

$  |f(z) | \le C e^ {2\pi M | z|  }$ \\

As in the classical case the passage is between Real and Complex . First  a  function F :$C^N$ $\longrightarrow$  C is   \underline{ entire } if each  coordinate function is. 
Let us define this notion  ( now  K is a convex body in $R^n$): An  entire function F ( as above) is   \underline{ "of exponential  type  K  "}   [17] :--

This is so if for every $\epsilon > 0 $ there  is a constant C =$C_{\epsilon}$ such that for every z in $C^N$   we have  
\begin{equation}\label{E:cong}
	|F(z) \le  C  e ^ {2\pi(1+\epsilon) ||z||_K } .  
\end{equation}   

For above Tnequality we recall   from Introduction  eq(3) \underline{ Minkowski gauge function }   and equivalent support function\\

Following [18 ] the class $\scriptsize E$ (K) is the class of all  functions of  exponential  type K\\

The  fact of the FT vanishing off a convex symmetric body  is explicit in this next Theorem from [18  ] It is   a main tool in  Sec2:\\

\underline{Theorem} [18]

Suppose F is in $L^2$ ( $R^n$ ).  Then F is the FT  of a  function vanishing off a symmetric convex body K iff  F is the restriction to $R^n$ of a function   in $\scriptsize E$ (K* )

\underline{Special case}  If   we      have an Interval  $ K=[-M. M]$    then it may be verified ( using the equivalence of the norm $||  y ||_K^{*}$   with the support function   $h_K(y)  $of K  noted in Sec 1.2) that  the inequality(5)  becomes  for y in R  called above as  Exponential Type M

\textbf{ Note:} In connection with   "Function of Exponential Type K"    above ,the authors of [6]  use the inequality (5) with  iy  where  y is in$R^n$ in place of  z in  $C^ n$  . \\

We  are now in a position to state and prove \underline{ Blaschke-Santalo Inequality }   Theorem 2.2  below from [6]

 We note the use below, of the polar K* of K in the  tool  or \underline{functional,$\rho$} used  there . We observe Condition (ii)  is  geared to  use the Theorem of  Paley-Wiener-Stein . Let us  define,from [6] this  basic tool,  or functional:--\\
 
\begin{equation}\label{E:conj}
\rho (K)  = inf ( \int_{R^n}   | F(x)|  ^2  dx  ) 
\end{equation}

the inf is over all " admissible"  continuous, $L^2$functions  F  : $R^n \longrightarrow$  C  such that

(i) $| F(0) |$  $\ge$  1 and (ii) the F T  ,  $\hat{F}$ ($\xi $)  = 0  for $\xi$  $\notin$ K  

Below we deal both the  Real Vector space $R^n$ and the complex one  $ C^n$;  in this section  B  denotes the closed Euclidean unit ball in $R^n$.\\

Also we denote in the next theorem   the measure or volume of  K in its linear span by "vol"\\
 
 \textbf{Theorem2.2 }[6]
 
 Let K  be  a convex body that is also balanced in $R^N$. Then:-
 
(a) P(K)   $\le$   P(B) , and (b) equality holds only if K is an ellipsoid.\\
  
  \underline{The proof} is done in these two key steps; the first is  easier
  
  \begin{equation}\label{E:conj}
     \rho (K) =  \dfrac{1}{vol K}
     \end{equation}
      and
    \begin{equation} \label{E:conj}
      \dfrac{\rho (B}{vol (B^*)}  \le\dfrac{\rho(K)}{vol (K^*)} 
      \end{equation}

  The Theorem  follows   at once from  combining  Eqs (4) and  ( 5). It is  in  proof of  eq(5) that   Stein's generalization of Paley- Wiener theorem is
  used. \\
  
  Let us prove  \underline{the easier,  eq(4)} .This follows from Plancherel's Theorem  (c) above , and Cauchy-Schwartz inequality .\\
  
  In fact,  it follows from  the Inverse FT  (b) , 
  
  F(0)  =$ \int_{K}  \hat{F}(\xi) d \xi $. 
  
  Hence  we may use condition(i) in def  of   $\rho(K) $  to write
  
  $|\int_{K} \hat{F} (\xi)  d \xi  |     ^2  \ge 1 $.
  
  Let us apply Cauchy Schwartz  to  see that left side of last   is
  $\le \int \chi _K d\xi  )(\int_K) | \hat{F} (\xi) |^2  d\xi $

  The first integral is  = vol(K)  and so  combined with Plancherel's Theorem  above(c)  again the last
  
   =  vol (K)  $\int_{R^n} | F(x)| ^2 dx $  
   
   and so  we have the last $\ge$1  from above. Thus,
   
$\int_{R^N}| F( x)| ^2 dx  \ge \frac{1}{ vol K}$
   
    As we may use " inf " over all  admissible F we get  
    
     $\rho (K)\ge\frac{1}{vol K}$. \\
     
      To    get  "  $\le$  let us consider  a specific admissible  F  namely:--  Let $\alpha$   be any complex number with  modulus  1 the    function  
       F( x)    = $\frac{\alpha}{ vol K} \int_{K} e ^{2\pi i x.\xi}  d\xi$ . Recallig Inverse  FT  above  we see that this function is just$\alpha / vol K $ times Inverse FT (b) of the function $1_K$  ie characteristic function of K.
       Hence we see that  the  F.T  of F   is    zero   off K  as required in condition (ii) and condition (i) $ |F  (0)|  \ge 1$   is  verified directly.  Finally we see that  $\int F(x)^2  =  1/vol K$  ; as F is admissible  we conclude  what we want ie$ \rho(K) \le  1/  vol (K) $ \\

     \underline{For eq( 5}) :  This is more involved than above: Let us  only  give the main ideas.\\
     A)  It  is enough to show that  
     \begin{equation}\label{E:conj}
     	\int |F(x)| ^2    dx   \ge  \rho(B)  \dfrac{ vol K^*}{vol B^* )}
    	\end{equation}
     For   taking " inf"  we get 
      $  \rho(K)  \ge  \rho(B)  \dfrac{ vol K^*}{vol B^* )}     $
      Here integral over $R^n$
      Remembering  value of $\rho(K)    $  from Eq(4)   gives the conclusion.
      
      Thus we must  prove above Ineq  (6) ;  we  see the Right Side next. For this write   K*   in polar coordinates

      \begin{equation}\label{E:conj}
       K^*= (t\theta : \theta \epsilon S^{n-1}  , 0\le  t \le ||\theta||_{K^*}  ^ {-1}  ) 
      \end{equation}  
      
      Hence fo;;ows  the integral in the next eq(8)   Here  let us recall  the dual  norm in the integral and $\omega_{n-1}$ is the vol of  $S^{N-1}  =   n Vol (B).$   This  Right Hand Side in (8) may be verified by  working out the  volume  of $K^*$from eq(7)in polar coordinates . 
      
      B) \begin{equation}\label{E:cong}
      	\dfrac{vol K^*}{vol B^*}   = \dfrac{1}{ \omega_{n-1} } \int_{S^{N-1}}||\theta||_{K*} ^{- n}  d\sigma(\theta),	
      \end{equation}

   we  next    combine above  two(6)  (7)   to get
   
  C) \begin{equation}\label{E:conj}
   \int |F(x)| ^2    dx   \ge  \rho(B)   \dfrac{1}{ \omega_{n-1} } \int_{S^{N-1}}||\theta||_{K*} ^{- n}  d\sigma(\theta)
\end{equation}

We  now   need to show this last  . It is here the authors  use the  Paley Wiener-Stein Theorem above,  Their ingenious idea is to find  function  called
 $ R_{\theta}(z) $,admissible for a multiple of B and  related to the initial F; in fact admissible  for  
 $   \rho (||\theta||_{K^*} B) $ .
  Let us  turn to this function  and  indicate the  main ideas  
\\
      
     \textbf{(1)}  Return to our F admissible for $\rho(K) $ . Without    any loss of generality we may assume F to be even . For  a $\theta  $ in $S^{n-1}$  we will  eventually  define our  function $R_{\theta}$  : $C^n  \longrightarrow C$ .\\
     
     . We continue with def. of  $R_{\theta }  $  :--
       Given  a  $\theta$  in $S^{n-1}$    define 
       $G_\theta : C\longrightarrow C$ by $G_\theta(z) = F(z\theta)$ .This function. we claim, is of exponential type $[-a, a]$  where $ a=||\theta||_K* ^{-1}$ . The reason is this:--
       The admissible  F ( domain $ R^n$) has the property :  its  F.T. vanishes off  K; hence it  satisfies the hypothesis  in the Paley- Wiener Stein Theorem above   .  We draw the conclusion : $  |F(\theta z)|  \le  e ^ {2\pi (1+\epsilon) ||z\theta|| _{K*}} $
       and the claim follows for  our   $G_\theta(z) = F(z\theta)$ 
       
     \textbf{ 2 )} Next $G_\theta (z)  =H_\theta (z ^2) $ for another function H due to  $G_\theta $   being even.
         At last we define  $ R_\theta : C^n \rightarrow C$ thus
         $R_\theta (z) $ =  $ H_\theta (\Sigma z_k ^2 )$;here  we write z =$ (z_k: 1\le k \le  n)$ \\
         
         We have
         
    \begin{equation}\label{E:conj}
     \int_ {R^n }|Fx)|^2  dx   =   \frac{1}{\omega_{n-1}} \int_{S^{n-1}}\int_{R^n}  |R_\theta(x) |^2   dx  d \sigma( \theta)   .
     \end{equation}

     The eq(10)   is  obtained thus:--
      We change to polar coordinates $(r,\theta) $ inthe integral in Left Side of 9(a) ; then replace as we may    $F(r\theta) $  by our newly defined function, $R_{\theta}
      $\\
       
      This enables us to   see the claim above about $R_\theta$ being admissible to the dilated   ball $||\theta||_K^*  B$   made above . and we have
      \begin{equation}\label{E:cong}
     \int_{R^n} |R_{\theta}(x)| ^2 dx  \ge\ ||\theta||_{K*}^{-n}\rho(B)
      \end{equation}

      \textbf{(3)} From this point on  the connection to  needed Eq(8)   via   eq(9) follows completing the proof QED
      
      The case of equality  part (b) is also non trivial; roughly it is reduced  to equality in eq(9) and an earlier characterization  whereby K    must be  an ellipsoid  QED.\\

       2.3. \underline{Another functional}
 This one called  $\eta (K)$ in [6] is similar to the earlier  "  $\rho (K)$" above but with differences  that make it harder.Let us  define it   below.
 
 \begin{equation}\label{E:conj}
  \eta(K)  = inf  \int_{R^n}  F(x) dx 
  \end{equation}
  
 The admissible functions F are continuous,    in $L_1$  and 
   
   (1) non negative  (2) F (0) $\ge$ 1   (3) the F.T $\tilde F$    vanishes off K
  
   There is  a \underline{conjecture} in [6]:
  
   $\eta(K) =\frac{2 ^n}{ vol K}$ \\
   
   The authors prove   this for  the cases of the Euclidean  unit Ball and the Cube. They remark, when the paper was in print, that this  was already proved by  Siegel  ( see references to this, and  above conjecture  in [6])\\
   
   Let  us indicate theorem from [6] originally   due  to Siegel for  these  bodies.The authors  remark that part of this is implicit in earlier work by Holt and Vaaler ( see reference to this  in [6] ) \\
   A tool here is \underline{ Poisson Summation Formula   [11].[17]  }
   
   Suppose f is in $L^1(R^d)$  .Then under  some  growth conditions   we have the equality:--
   
   \begin{equation}\label{E:conj}
   \Sigma f(m)    =  \Sigma \hat f(m) 
   \end{equation}
   
   the sum over the lattice   $  Z^d  $ 
   
   Following  [20]  in another context  we refer to the 2 sides of this Formu;a by " Function Side"  and "Fourier Side"\\ 
   
   We  use also function  f  called \underline{" sinc squared function"}
   
   namely   $f(x)= [\frac{sin \pi x }{\pi  x} ] ^2$   and its integral over R  is 1
   Its FT  is the  triangle function $ g(x)= 1- |x|  $   for  $ |x| \le 1$  and zero  elsewhere\\

  \underline{2.4 Theorem  }   
   
   In case of K = $ B_2$  or $B_\infty $ then $  \eta(K) =\frac{2 ^n}{ vol K}$ 
   
   Let us observe that the   case of the cube Q =  $B_\infty$ is easier   and  we  only consider this  here\\ 
   
  \underline{ Proof:}
   
    So, let K be the\underline{ Unit   Cube  Q}   .
  
    Since  vol (Q)  =$2^n$  we have to show that 
    $\eta Q =1$

      Firstly,starting with an admissible function   F for $\eta Q$, we use the Poisson Summation Formula Eq(15)  recalled  above.
    
   The  ( Left ,or) Function side" " has  all non  negative terms  ( due to $F(x) \ge 0$ ), and so this sum is $\ge F(0) \ge1 $   by assumption on F. In the ( equal   )" Fourier side " only $\hat F(0) $ survives. rest are zero due to$ \hat F$    vanishing off  interior of K =  Cube.  We  have thus the outcome:--
    
    $ 1\le \hat F (0)$   = $ \int_{R^n}  F( x)) dx $  and so  taking "inf"  over all admissible F  gives  $\eta(Q)\ge 1$. \\
    
    To get the opposite inequality we consider a  function
    
     F(x)=  $\Pi (\frac {sin \pi x_n } {\pi  x_n}  )^2  $,
    product over $n=1  to N$  .  Due to the  F.T .of the "  Sinc Squared" function recalled above ,this function  has support  contained in   Q  whence F is admissible . Now,  integrating it gives  1  ( using   the fact that integral of    each factor  of this  =1  as  observed above)
     Thus we have  $ \eta Q  \le 1$ as required . Hence this part thereby completing the proof  QED .\\

   \section{The lower bound for Zonoids [9]}	

   Now let us  consider the result  from [9] proving the   lower bound  in the class of  zonoids. The paper is " essentially   self contained"(ie.using  standard  methods and a result on the volume of a zonoid   from[17]) . There are 3 Lemmas  ; however the proofs are non trivial.

    3.1 Let us  record a few tools needed below:--
   
    \underline{ Brunn-Minkowski Inequality}    [ 4 ]   [8  ] [17] [20]   thatwe recall  at Lemma 3.5   below where its needed
    
     Let us  next  recall   the term, \underline{zonoid}.  From  [7]  for  ex. this is  due to the 1940 Liapunov's Theorem  :--   It is a convex and compact  set A.  It  is  the range of a non atomic vector measure into $ R^n$; or,the closed convex hull of the range  for general , nor necessarily non atomic, ( vector) measures. It is  proved in [7] that there is then a (unique)  positive, even finite measure$\mu$ on $S^{n-1}$  that generates A; roughly,  A  is   the set of all " "averages" taken with respect to $\mu$, or, "moments" of the identity function  on $S^{n-1}$. Further we have  from[ 7] for the dual norm above,or   support function (of a balanced zonoid)  that (with this $\mu $) 
     
     \begin{equation}\label{E:cong}
      ||y||_{A^*}  = \frac{1}{2}  \int_{S^{n-1}} |<x, y>|  d \mu (x)
     \end{equation}	  
     	 
     	    holds for   y$\epsilon$  $R^n$  .  This is the  modern  definition of "Zonoid"[9]  ,[10],[19]   and  used  below.

    Also   let us  note here the \underline{ projection  $P_x^\perp$ }  ; given x in $R^n$  let H(x) be the hyperplane thru O and  perpendicular to vector x, This projection map is from $R^n$  onto H(x).  
    We  give examples to illustrate the Lemmas  and hope these are of interest.
    \\
    
    We can now come to   main   result  :\\

\textbf{3.2Theorem}([9])   Let A  be any  zonoid  in $R^n$;  then P(A) $\ge$ $ 4^n$  /$n!$
with equality  (if and) only if  A is an  n cube.\\

 The  proof is by induction on "n" ; the 3 Lemmas   take care of the passage  from "n-1" to "n".  We defined above  the  supporting measure  $\mu$  and \underline{projection } $P_x^\perp$    .Also the fact that the Projection of the  zonoid  A is  also a zonoid of lower dimension ;which facilitates the Induction  in  the proof  . The idea is this: Given zonoid  A    we find  another , $A_1$   of  dimension (n-1) , The  crux is to prove :--
 
 $ P(A) \ge  \dfrac{4}{n} P( A_1)$   and Induction applies  to $A_1$\\

  \underline{Lemma 3.3} ([9] Lemma1) Let A be a  zonoid in $R^n$  with supporting measure $\mu$. Then 
  
  \begin{equation}\label{E:cong}
 (n+1)|A|\int_{S^ {n-1}}[ \int_{A ^*} | <x, y>| dy ]  d\mu (x)    
 	= 2 |A ^ * | \int_ {S^{n-1}} | P_x^\perp   (A)| d\mu (x)
 \end{equation}	
 	 
 	 In particular, for some $x_0$ in $S_{n-1}$ we  have
 	 
 	 \begin{equation}\label{E:cong}
 	 (n+1) |A|  \int_{A ^*} |( x_0 , y) | dy \ge   2|A^* || P_{x_0} A| . 
 	 \end{equation}
 	  \\
 	 
 	\underline{ Proof:}
 	 
 	  As mentioned before the proof uses for the formula from [15] for volume of zonoids, namely  
 	 
 	 $|A| =n^ {-1} \int_{ S^{n-1}} | P_x ^\perp A| d\mu (x) $
 	 
 	 We use the formula   for the zonoid  A  for $||y||_A* $ above Eq(14) .\\
 	 
 	 Let us start with the  integral in Left Hand Side in eq(13):
 	 
 	 $\int_ {S^{n-1}} [\int_{A ^*}  |<x, y>|  dy ] d\mu ( x) $   this gives  by  Fubini's Theorem and the formula in (Sec 1.2) for  $||y||_ A^*$ that this  is in turn
 	 =  $ 2 \int _{A^*}   ||y|_ {A ^ * }  dy  $
 	 
 	 =  $ 2 \int_{A^*} [\int_{0}^{ ||y||_ {A ^*}}  dt  ] dy $ (  we note the trick of writing  x as  $ \int_ {0}^x  1 dt $). Also
 	  we may rewrite this last (again by  Fubini's Theorem) 
 	 
 	 =$ 2\int_{0}^{1} | (y\epsilon  R^n  : t\le  ||y||_{A^*} \le 1) |  dt$\\
 	
 	 Let  us  note the integrand as between  volumes  of the    sets 
 	  $ A^*$  and  $ t A^*$  ; this gives   $ |A^*| ( 1- t^n)$ .  integrating this  last  we get  
 	 
 	 $ \frac{2n}{ n +1} | A^ *| $.  With this in hand  and the formula  above   for the volume of zonoid  A   we get the Right Side . Hence the Lemma
 	 
 	 The " in particular  " follows from standard  arguments
 	 \\

 	 The next Lemma is adapted from    earlier    authors ' paper on Inequalities and is also    self contained. We note the concavity hypothesis .
 	
    It  uses the function   
 	$  (x-y)_+ $  this is =  x-y  if $x \ge y$  and =0  otherwise.\\

\underline{Lemma3.4 [9] Lemma2}  The   conditions on
 	 f:  $R_{+} \longrightarrow R_{+}$   are:

  (i)  f(0) =1 f not identically zero  and that(ii)   $ f^{1/p}$   is concave  on the support set of f where $p>0$.   Then we have the inequality, 
 
 \begin{equation} \label{E:cong}
    \int_{0}^{\infty} t f(t) dt \le  \dfrac{p+1}{p+2} [ \int_{0}^{\infty}  f(t) dt  ]   ^ 2 .
  \end{equation} 
  Equality holds    iff  f(t)  =  $(1-at)_{+}$ $ ^p $  for some a  $>$ 0   \\
 
   \underline{Proof:}  We begin by defining a  positive number  a  by  equating  following  first and last:
 $\int _{0 }^{\infty} f(t)  dt =\int_{0}^{\infty} (1-at)_+^pdt   $  and (this last)=  $\frac{1}{ a(p + 1)} $
 
 Next let  g(x)  = $ f(x)- ( 1-ax)_+ ^ p$ Then  we have the following   for   it remembering the assumptions on  f:  
 
  g(0)= 0,  $\int_{0}^{\infty} g(x) dx$  =0  ( due to the equal integrals above)  .   Also the  concavity  of  the function  $ f^{1/p}$has the consequence 
 
 \begin{equation} \label{E:conj}
   \int_{x} ^{\infty}  g(t) dt  \le 0   for all        x\ge 0 
   \end{equation} 
  
  Now we may compute the integral on left  hand side of the statement of this Lemma:
  
  $\int_{0}^{\infty} t f(t) dt  $ 
  
   =  $	\int_{0}^{\infty}[\int_{x}^{\infty} f(t) dt ] dx $

   $\le \int_{0}^{\infty} [ \int_{x}^{\infty} (1- at)_+ ^p dt] dx $   ( we used eq (14)).
   
   Now we may compute this last  to be
   
   = $\frac{1}{(p+1)(p+2)a ^2 } $
   
   =$\frac{p +1}{p+ 2}[ \int_{0}^{\infty}    f(t) dt ] ^2$  ( recalling the equation at the start of the proof that defines a  or rather 1/a)
   
    Hence the inequality in this Lemma.
    
    As for equality  we look at  the  two( double) integrals  above  with $\le"$  . Therefore  this equality   forces  that of the inner integrals  for every x  in  $[0, \infty) $. This last  is only possible  with equality of the integrands.  Hence the Lemma .\\
    
     For the next Lemma  use   a  consequence of the \underline{Brunn-Minkowski Theorem} ([3], [4] ,[8]  , [13]  [20])    that we  recall it  here 
     
       It is used in the form:--
     Let B be a symmetric convex body , x  in $S^{n-1}$  , recall   hyperplane  $H_t= ( y \epsilon  B : (x,y)= t )$   and  volume of slice  $  g(t)  =  |  B\cap H_t| $.  Then the function
     $ g^ {1/ n-1}  $  is  concave on its  domain $[- a, a]$   where
     $a= ||x||_ B  ^{*}$
    
    For our needs it  means that the function  g (=measure of  slice of body B) defined below  in Proof of Lemma 3.5 has the property :--
    	 $g ^ {1/ n-1}  $  is concave  on its support   \\

	\underline{Lemma 3.5}([9]Lemma3)
 	
 Let B be convex body , symmetric about the Origin. For  x $\epsilon$  $S_{n-1}$
 	let B(x)=  (y $\epsilon$ B : $ < x, y>$  =0)  .Then
 	
 \begin{equation} \label{E:cong}
 	\int_{B} |<x, y>|  dy \le \dfrac{n}{2(n+1)}\dfrac {|B|^2}{|B(x)|}
 	\end{equation}
 	
 	The conditions on equality   (  with the x in statement) are  iff   there is a $y \epsilon R^n$  so that  B=  conv (y, -y,   B(x)).\\
 	
      \underline{Proof:} Fix x ; then  and  let t $\epsilon$ R ,and  g(t)=$| ( y \epsilon  B: <x, y  =t) | $
     
     Now g is an even function, g(0) =B( x)   and gis supported on the Interval,  I = $[- ||x||_  B^*,   ||x||_ B^*] $
     
     We  may thus  apply  Brunn's    Theorem  and obtain concavity of $ g^\frac {1}{n-1}$  in the interval  I .
     
       This allows us to  apply    Lemma3.4  below.
    
      Now   we have  on the one hand 
      
 $ |B|  =  2\int_{0}^{\infty} g(t) dt$   and on the other
       
       $ \int_{B} | < x, y >|  dy  =  2 \int _{0}^{\infty} t g(t) dt$
       ( this last by rewriting   the left integral using  the  trick    $ x= \int _{0}^{x} 1  dt  $    we used in proof of Lemma 3.3   with Fubini's Theorem )
       
       The result is that we are able to apply Lemma 3.4   with  the function  "  $f= g/g(o) $and  $ p= n-1 $ to get the inequality.
       In    detail, the Left side  in Lemma3.4 is   now
       
       $\int_{0}^{\infty}  t g(t)/ g(0) dt$ and the Right side is   
       
     $\frac{n}{ n+ 1}   [ \int _{0}^{\infty}  g(t)/ g(0)  dt]^2 $
       
       Remembering the equations ("On the One hand...)   above   we  get the desired inequality.  
       
       As  for equality we may apply this case  from Lemma 3.4 for our f = g/ g(0) and the constant a there 
        =$ 1/||x||_B^*$. Thanks to the Hahn Banach Theorem there is a y  in B   such that this norm  $||x||_B^*  =   | <y, x> |$. Now let $B_1  $= convex hull of  $(  y, -y , B(x))$ ; then we have$ B_1 \subset B$  and we claim that these are equal  which will finish the "only if" part. The claim follows from the following chain  of inequalities:
      
       $|B| \ge  |B_1|  = \frac{2}{n } ||x||_B^* =2\int_{0}^{\infty}  g(t) dt  =  | B|  $  which compels equality of B  and  $B_1.$
       
       The "if "part may also be  verified\\
       
         \underline{ Remark Ex 5}  Let us  note   in Lemma 3.5 above that the equality condition is for special  x in $S^{n-1}$.  Also we  note the  set in "Equality " part  writes    B  as   double  cone"   B = $conv ( B(x),  y, -y)$  This last is usually not  a zonoid [7]:--
         
         In exs  below  let us  look at these :--
         
         A)  The inequality  (18)  and if  equal?
         
         B)   Find the Polar set 
         
         c)           Verify  Volume Product  Inequality esp  Lower bound and also if Equal?

         \underline {Ex5a} Let B=$ B_2$   the Euclidean unit disc and x=$( cos\alpha, s in\alpha)$  in $S^1$. Then there is strict Inequality in(15) as  B  cannot be a cone of this , or any sort.
         
         In fact, we have  B(x)=  a diameter  hence of  length2  , and $| B| =\pi $ ;this gives
         
         Right Hand Side(RHS)= $  \frac{2}{6} \frac{\pi^2}{  2}$=$\pi^2 /6$  whereas
         
          for   LHS   with x as  above let us consider
           $y=r ( cos\theta, sin\theta)$  and so 
          
          $<x, y>$ =$r(cos\theta cos\alpha   +  sin\theta sin\alpha  )$  =  $  r cos(\theta -\alpha) ) $ Hence the integral in LHS  =
          $\int_B r | cos(\theta- \alpha)| r drd\theta$
          
           where we used polar coordinates for the double integral  This works out to $4/3  < \pi^2  /6 $  =RHS\\
          
          \underline { Ex 5b}  Let  B =  conv  $ (- e_1,  e_1,  e_2, -e_2  )$   with  x=  $ e_2 $   Then we have  equality in (15). The  set B is now the "diamond" or the  unit closed ball in norm of $l_1$ ie the set
       $(  (x, y)  :  |x| +  |y| \le  1)$. 
       
          We have $ | B  |$ =2  and also $| B(x)| $ =2  , since
          B(x) is again the segment   with length 2    and so 
          
          RHS=  $  (2/6) .2^2  /2$  =  2/3.
          
           As for LHS  writing  $x=(0, 1 )=e_2$ and $y= (y_1, y_2) $ we  may calculate $ <x, y>$ =$ y_2$ .Hence the  LHS integral is $\int_{B} | y_2|  dy$.
          This time we use Cartesian coordinates and get  2/3 .
          
          so = RHS.
          
          With 
          \\
        
  \underline{ Ex 5c} In $R^3$ let $B(x)$  be the (closed )unit Euclidean disk in the plane  xy  ie  in $x_1   x_2$ plane  and  $ y= e_3  = ( 0, 0, 1)$ . Now the   body C=   conv $ ( B(x),  y, -y)$  is a  double  cone  and  is not a zonoid[7] as  it is not decomposable   Let us note   we need the   defining  inequality
  $ \sqrt  {x_1^2  + x_2^2 }   + \ | x_3| \le 1$
  
   (A)Ineq(18) is satisfied  with equality . both sides $\pi/6$   As usual Right Hand Side is easier  .The Left  needs the    calculation: ( multiple integral)  $ 2\int_{z\ge 0}  z  d V$  where we used the usual ( x, y , z)  and  note the defining inequality for C
   
   (B) Polar is a cylinder with height 2  and  unit disk  base. Namely
   $ y_1^2  +y_2^2   \le1   $    and $|y_3 \le1$
   
   (C)Volume Product 
  \\
         
         The analog   of Ex 5b in $R^3$  is  the Octahedron,   the Unit Ball   $B_1$  in norm of $l_1$; this too is not a zonoid as we noted before Sec 1.2  but  we may verify likewise  for $ x= e_3 $ that  Ineq(15)  with equality.\\
         
         \underline{	ex5d}   For another  example with equality  in  $R^3$  let  D be  the unit    square again in xy plane ie  in   $ x _{1}  ,  x_{2}$   plane     and $C=conv(\pm e_3,  D) $   a sort of "double   "pyramid ".  Then ( according to Lemma 3,5)  the polar $ C^0$   is   a  Cube   . 
         
         B) However we  may check it is  the  the Cartesin   Product of the Disk  $B_1$ with the interval $[e_3,  - e_3]$.

         A) We may verify that equality holds; both sides equal$2/3$ The righr hand side is easier . The left by a calculation of multiole integral,  $\int _C   y_2  d  V) $  .
         
         \underline{ Remark 5d} We emphasize the  two examples in Equality case of Lemma 3.5  are non zonoids  due to triangular    faces but the polars are.
           \\ 
       
       \underline{Proof of  the Theorem} 3.2. This is by induction . If n=1  the result is true   easily as we see from $A=[ -a, a] $ and so$ A^* $=  $[- a^{-1}, a ^ {-1}] $  and P(A) =$2a.2a^{-1} =4  $ . So assume that the Theorem is true for  n-1. To  go from this to n   is done  by combining Lemma3.3 and 3.5  to    get the following:
      
      $(n+1 |A| \int _{A^*}  | < x, y>|  dy )   \le \frac{n|A^*|^2}{2 | A^*( x)| } $ by Lemma 3.5  and also
   
     $2|P_x ^\perp  A  |A^*|  \le  (n+1 ) |A| \int _{A^*}  | < x, y>|  dy $   by  Lemma 3.3
    
    where  the " x" is one given by the    "in particular "  Eq(18) in  Lemma 3.3.  Putting these last two together ( and canceling $|A^{*}$) we have
    
    $ 2  | P_x ^{\perp} (A) | \le  \dfrac{n |A| |A^{*}|} {2|A^ *(x) |}$
    
    Let us note the fact  (we may verify):--$P_x^\perp A =(A^*(x))^* $
    
     With its help and recalling $P(A) = |A| |A^*  |$  the     outcome   of  the last  inequality is this:
     
      $ P(A) \ge (4/n) P(P_x ^\perp A )$.\\
    Now   we may apply induction hypothesis to Projection of the zonoid   A , namely $P_x ^\perp   A $ is a zonoid ,  and has dimension n-1  so  that  we  have, thus
    
    $P(A) \ge (4/n) (4 ^ {n-1}) /  ( n-1)! $ = $\frac{4^n}{ n!}$  completing the induction and the proof of Theorem 3.2    except  for the case of Equality.
    
    This too is by Induction   similarly to above, and hence the Theorem QED

 \section{Perturbations  of unconditional bodies  [ 12] }
   
    In this Section ,the idea is  "to perturb  slightly , whatever that means"  in [ 18 ] and  used effectively in [12  ].  
    
    We prove only one  result   Theorem4.2 from[12]    below.

    The   tool used for the perturbation is the   classical  \underline{Banach-Mazur distance} , now  between two convex bodies K and L  ( here   GL(n) is the space  of all isomorphism  of $R^n$):-
   
    \begin{equation}\label{E:cong}
    d_{BM}  (K, L)=  inf ( d\ge 1 :L\subset TK \subset dL \exists  T \epsilon GL(n) )
    \end{equation}
    
     The above   does not define a metric ; if  K= L  we have  $d_{BM}  (K, L)  $=1 .In words of [12]," it is multiplicative  not ( sub) additive ". To get a metric we  need  log d .   Another well known metric ,also used in [12] on the set of ( nonempty ) compact sets of $R^n$ is the Hausdorff metric. However  we  use the  above.  Since this isthe only  distance used here we  simply write  " d  (K,L)  " \\

     We need concept of\underline{ unconditional body} K .It goes thus:--
     
     These   K   have the property:
     $y = ( y_1, y_2,........y_n )   \epsilon $  K  whenever 
     
     $x = ( x_1, x_2, x_...... x_n) \epsilon$  K   where each $y_i$  is   $\pm x_i$ .  As  explained in [11], this  means that K is stable under  reflections of all coordinate hyper planes, $x_i$  =0
     We use  in  Sec 5  the  corresponding notion of Unconditional function

     We  need next   the concept of \underline{Hanner Polytope }[12]:This is a convex body  H  that  is 1-dimensional, or it is the $l_1$  or $l_{\infty}$ sum of two lower dimensional Hanner Polytopes. We recall these last two sums  for  convex sets  A  and  B .
     Their ${l_1}$  sum is the convex hull of $A \cup B$  and the ${l_{\infty}}$ sum is  A  + B ie. Minkowski  sum \\

    The\underline{ object}   of  [12]  is to prove  \underline{" stability  of the lower bound "}  in   the class of    symmetric convex bodies K ,in the sense ([12]Theorem2)  namely
    
    \begin{equation}\label{E:cong}
   P(K) \ge   ( 1 + c(n)\epsilon) P(B^n_ {\infty})
   \end{equation}
    holds, provided K itself is sufficiently close to an unconditional body for  appropriate constants depending on n (Theorem4.2  below) and of course, $\epsilon  >0$ and arbitrarily small. We only give  indication of this main result
   
    To avoid repetition let us consider bodies such as K, L, H, $H_0$  to  be in $R^n$  unless otherwise stated explicitly.\\

  Very roughly , the following  Theorem is a sort of "denseness" of Hanner polytopes given the condition about  K being " close  to" the $B_n^{\infty}$  but in sense of  its volume product\\, 

\underline{4.1Theorem}  [12, Theorem  1 ]

Let  K be  unconditional .

If    for some  
\begin{equation}\label{E:cong}
   0 < \epsilon\le\epsilon (n) ,  we  have  
| P ( K)  - P (B^n _{\infty} ) |\le\epsilon
\end{equation}  

 then there is a  Hanner polytope H such that
  
 \begin{equation} \label{E:cong}
 d ( K, H) \le   1 +  c (n)  \epsilon.
  \end{equation} 
  
  The constants  $\epsilon$(n ) , c(n) $>$0   depend  only on n.
  The authors remark that here $ n\ge3$   as   this Theorem was proved earlier  in case n=2 .

The proof  is non trivial  and depends on  a few steps; Lemmas 1  and 2. The  Hausdorff distance  comes back in an equivalent formulation of Theorem 4.1   
 We already noted  the  essential use of the  Hanner Polytopes .  
 
 For next result of[12]  recall that $ x^\perp $  is the Hyperplane  thru O and perpendicular to x  so that $ K\cap  x ^\perp$ is a slice of K  .Hence   the conclusion compares  this slice of K  with  the (n-1)  ball on $|| ||_\infty $  norm  . Both are in sense of volume product.

 \underline{4.1A  Lemma }[ 12  Lemma1]

 Let $\epsilon >0$ and K unconditional such that
 
$P ( K)  \le  ( 1 +\epsilon) P( B _1) $  then

$P( K \cap e_j ^\perp) \le ( 1 + n \epsilon ) P(B_1^{n-1}   )  $    for every  j  between 1 and n  ( now  $ n \ge2$)

The proof is  rather involved  .\\

	 In   the next Lemma,we note  the  hypothesis ,  equality of  certain  sections of K   equal to  those of the Ball in $B_{\infty} ^ n$  that we write as $ B_{ \infty}$ , and likewise  $B_1$ and dropping of    " unconditional" .   By   the supposition on K and on p    we must  necessarily   have $    t\le1 $.  We  also  note  volume products are  compared  in conclusion .\\
	
	\underline{4.1B Lemma}[12, Lemma2]
	
	Suppose that $K \subset B_{\infty} $  with $n\ge3$  and
	
	that $K \cap e_j^\perp =  B_{\infty} \cap  e_j ^ \perp$
	
	holds for every j between 1 and n   and  that
	 $ p= t( 1, 1, 1,   ..1)$   is in the boundary of  K where $t\ge 0$ .  Then we have  the consequence:
	 
	 $|K| | K^*|  \ge  ( 1 + c . (1-t)) |B_1| | B_{\infty} $.
	 Here  c = c(n)  is positive and depends only  on n

The\underline{ stability of  the  lower bound} for   convex  bodies under consideration  is achieved in the following sense;  as we see this takes some effort--\\

\underline{4.2 Theorem}  [12, Theorem2]

Let K  have center of symmetry O ,and satisfy the following:--

  $ min( d  (K, L) : L $ is an unconditional convex body) $=1 +\epsilon$   
  
  where $ 0 <  \epsilon   <  \epsilon$ (n). Then
Inequality (22)  holds  ;namely 

$P(K) \ge (1+c(n)\epsilon)  P(B^n _\infty) $. here $c(n)$  depends only on  n

 \underline{Proof}:  this  too, involves a number  of steps, especially the next Theorem 3  where the hypothesis of above Theorem is now with a Hanner polytope in place of L, namely   the following ( proved earlier by the    1st  author ) that we quote without proof:\\

\underline{4.3Theorem( [12] Theorem 3)}   Let K be a symmetric  body that is close to some Hanner Polytope  namely:--

  min  $(d(K, H) $ :  H is a Hanner Polytope  )  = $ 1 +\epsilon $ , for some $0< \epsilon \le  \epsilon (n)$

 then
  $P(K) \ge P( B_1^n)  + \epsilon c(n)$   for   some $ c(n)$  depends  on  n only\\
  
  First let$\gamma_n$  be so chosen that( with $\epsilon$)  both Theorem 1 and 3 are satisfied.   ie  we  have the following   ( roughly contrary cases).In both $\epsilon \le  \gamma_n$  :--
 
  Let us  use the phrase" $1+\epsilon$   close " to mean in sense of B-M distance,
 
 (a) If  K is $1+\epsilon$  close to\underline{ some} Hanner Polytope  then  
  
  $  P(K) \ge (  1+ \alpha_n  \epsilon) P( B^n_\infty)  $
  
  (b)  If on the other hand  our unconditional body K is  " $1+\epsilon$   far away "     from any Hanner polytope then the  above Inequality in case(a) is satisfied  but with $\beta_n$    instead of $\alpha_n$  there ; 
   $ \alpha_n   $and   $\beta_n $  only depend on n.
  
  \underline{Proof of Theorem  4. 2} depends on  both Theorems4.1  and4.3  but is not immediate.   Roughly, we again have two cases  in which  the above " contrary cases"  may be applied.  Let us only give a sketch.

   Given our K  we   first choose  unconditional L   that is $1+\epsilon$  close to K then choose  a Hanner polytope H  that is $1+\delta$    close to L  ; these are to be smallest possible   choices .
  
    There are  again two cases to consider; they depend on  numbers  $\delta$and $ \gamma_n$ . These stem from Th 4.3   and 4.1 and  correspond to this \underline{dichotomy} on   : case a) $\delta \le 2\gamma_n / 3$   or   caseb) not      :--\\ 
  
  Either there exists a Hanner   Polytope  H that is$1+ \epsilon$ close to K   or  there is no such  H  . Accordingly we choose  $\gamma_n  >0$  dependent on $\epsilon$ to satisfy both Th 4.1  and 4.3   Namely
  the constant$\gamma_n$ satisfies these statements:
  
  \textbf{(c)} If K is a symmetric body  with    $d ( K, H)   =  1+\epsilon$    (with $\epsilon \le \gamma_n $ )  
  
   then:   $ P(K) \ge (1 + \alpha_n  \epsilon) P( B_1) and  \alpha_n  $  depends only on n 
   
    \textbf{(d)}  is opposite to(c)   ie  if K is an unconditional body with the property : for every  Hanner Polytope we have
   
     $ d _{BM} ( K, H) \ge 1+\epsilon  , with ( \epsilon  \le  \gamma_n) $  \textbf{b1}  
     
     then:  $ P(K) \ge ( 1 +\beta_n \epsilon) P(B_1)$
     
     $\beta_n$ depends only  on n\\

  The choice of $\gamma_n$  is called  by  the authors of [12]   " a  threshold  "of $\epsilon$  "
  
 To prove Theorem 4.2 we start with our  convex and symmetric  body K  ; then  choose   $\epsilon >0$ and (" smallest  ") unconditional body  L  so that    $d_BM  (K, L)=  1 +\epsilon$ .  Further we may choose  ( again smallest  among its class ) Hanner Polytope $ H_0$ , and  now  $\delta >0$   so that  $d_BM( L, H_0)  = 1 +\delta$.
  
   Again there are  2 cases:  (i) $\delta \le  \dfrac{\gamma_n}{3}  $   and case(ii) the opposite . They  fit  quite well into the dichotomy  above  \textbf{(a)    ( b)}  .   Due to the different  conditions the  way to handle is also different , Case(ii)   more involved than Case(i)\\
  
 \underline{ case(i)}  Assume $\delta \le  \gamma_n /3$  so that we have,
 $d (L, H_0)  \le 1 +  \gamma_n /3$. Then we have these inequalities:  First we use the" multiplicative"  nature of the Banach Mazur Distance  to pass from  K to  L to $H_0 $   to get
 
 $d( K, H_0)  \le ( 1+\epsilon) (1 +\delta)$
 
 $\le  1 +\epsilon +( 1 + \gamma_n  /) $ using our assumption
 
 $\le  1 +\gamma_n$  as we have chosen as we may $\epsilon_n  \le  min ( 1,   \gamma_n  /3)$
 
 Hence the condition in \textbf{( c )} is satisfied and so  we may conclude Equation there  :that  $ P(K) \ge  ( 1  +\alpha_n  \epsilon)   P(B_1)  $
 This is  what is wanted  \\
 
 \underline{case(ii)} As noted above this is handled differently:
 So assume that$\delta  >\gamma_n /3 $ and so $d_{BM} ( K, H)
 >  1+ \gamma_n /3$   for  every Hanner Polyope  H.  We  will  require
 $ \epsilon_n  \le  min( \beta_n \gamma_n,   1/2n ) $.
   Rcall choice of   L from the  start  of this proof.   As   there is some T in GL(n) so that  the inclusion 
  
  $ L \subset T(K) \subset( 1 +  \epsilon ) L$)holds,  we  have
 
 $P(K)  = P(TK)  \ge |L|   | (1 +\epsilon)^{-1} L | $
 
 this last is   $({ 1 + \epsilon)^{-n} }  P(L) $  due to rule for${ c B}$  for  any body  B ,	which in turn  gives
  $\ge (1- \epsilon)^n  P(L) $  and this is
 	
 	$\ge ( 1-n\epsilon)  P(L)$  thanks to a  limited Binomial theorem estimate.
 	
 	 Moreover we are entitled to use \textbf{( d )} above to L, therefore  we   have
 	 
 	 $P( L)  \ge  ( 1 + \beta_n \gamma_n  /3) P( B_1)$
 	 Now  we may put these  last 2 inequalities  together  and conclude after some calculations  that  
 	 
 	 $ P(K)  \ge ( 1  +n\epsilon) P( B_1)$
 	 
 	 To end let write   $\epsilon_n =min (\gamma_n/3,  \beta_n \gamma_n /12n, 1/2n)$  and
 	 
 	 $\tau_n $ = min($\alpha_n, n $)  and we have
 	 
 	 $P( K )  \ge  ( 1  +\tau_n \epsilon)P( B_1)$ provided $\epsilon  \le  \epsilon_n$
 	 This completes the proof of Theorem 4.2    QED
 
 \section{ Generalizations; functional version  [ 10]}

  The results in this Section may be   compared  with the  earlier Prekopa-Leindler  Inequality[4 ], [ 8 ], [16] [ 18 ] that generalizes   Brunn-Minkowski 's.   These  new generalizations ( as indicated in Section title )  all depend on the notion of  "\underline{Polar of a function" } ;  see   the references   [1]  [2]  [3]  [10]  also[15].     We content ourselves  in discussing  an Inequality  called in [10] as   " Ball's Inequality"( see  Eq(27 )).  Since this last   notion of "Polar" is at first not intuitive  (to us ) let us   discuss it from[1] [2][3]   then  go to Ball's Inequality.  We   end by briefly  discussing two Inequalities:- the forst EQ(31) from{1}  [2]  when specialized  gives   Blaschke Santalo Inequality  . The second eq(32)  from [10]  likewise implies Balls' .

    Given an even, non negative valued (measurable)   function f: $R^n$ $\longrightarrow$  $R_ +$  its \underline{polar  $f^0$} is   defined  thus:-
  
  \begin{equation} \label{E:cong}
  	f^0 (x)  = inf (\dfrac{e^ {-<x,y>}}{f(y)}  : y\epsilon R^n  )
  \end{equation}
  
  This  definition was initially given for   Log Concave Functions [1], [2], [3]   and later as above in [10];see also  the classic[15] We  we  note from  [1],[2] [3] ,[10]  the relevance of this definition to classical Legendre Transform. \\ 
  
    Let us  mention that the  references [1] ,[2] explain how this  definition  meets  \textbf{criteria  for " polar"}  . Let us  give these in  a chart:\\

    	Let us denote the euclidean unit ball by$B_2$ and the Gaussian function by $g_2$.  All these in the chart can be found in  [1]  [2]
    
    \begin{tabular}{|l| r|r| r|r  |r |r |r|r |r |} \hline
    	Name  &   $  $   &   &   &   &  &  &    \\  \hline
    	Set   & A & $\subset$ &B &$ \longrightarrow$ & $ B  ^ o $ & $\subset$ & $ A  ^o $	\\	\hline
    	function &f   & $\le$ &g & $\longrightarrow$  & $ g^o$ & $\le$ & $f^o$ \\  \hline
    	both  & $A^{oo}$ & = & A &  & $ f^{oo}$ & =& f\\ \hline
    	Selfpolar&$B_2$  &$  $  & $    $ &$   $  & $g_2$  &  & \\ \hline
    \end{tabular}
    \\

 Also ( see [2] ) it is shown that this  definition of Polar  of a function is  essentially the  one   for a map  on  a  suitable space  of  functions that  preserves  above   conditions.  let us  note some \textbf{Facts}   from [1], [3]  about the Polar of  functions on $R^n$:-- \\
  
  Perhaps we may remark  in above connection, that    the polar of the  Characteristic function  of K is,   as with any function is also a function  but  " not   another of same kind".
  For the first Fact let us  recall the "dual norm "  in sec2  recalled for ease of reading  $ ||y||_{K^*}    =  sup  (  |( x, y) : x\epsilon K   ) $ \\ 
  
  \underline{ Fact1} Let   K be a convex symmetric body  and  f= $\chi_K$   its characteristic function ( ie  =1  on K and  0   off it), Then its  polar  $f^0$ is  $ e ^  {- ||x||_{K ^*} } $ 
  
  We  may see  this ( from  eq(24) , noting first that y  there must be in K . Next  that the minus  in numerator sends that exponential to denominator , therefore we must take " sup".    Now  with this  ( applied to   $< x, y>$   with  y   varying in K )   we see  the appearance in exponential, of " dual norm"  from   Intro  as stated.   To be sure we may also verify  on similar lines   the "involution" in Chart  namely:  we have  $f^{00}0=f  $.  \\
  
  \underline{Fact2} Let    K  be as before  and with its polar set K*  let us   denote  for brevity their  corresponding norms  by $ ||    .||$  and  $ || .||^*$   . Now let   f ( x) = $e ^{-  ||x ||^2 /2} $   then its polar is  $f^0  (x)$  = 
  $ e ^ {- ||x|| ^ {* 2}  /2 } $ \\
  
  This too may be verified .\\
  
  Let us  proceed now to the Inequality    called in  [10 ], as  \underline{"Ball's  Inequality"}; it is stronger   than the classical one above  eq(1)  for  unconditional bodies.\\
  
 \underline{Theorem 5.1} Let  K   be such a body and $\omega $   the volume of the n-dimensional  Euclidean Ball. Then
  
  \begin{equation} \label{E:cong}
  	\int_{K} \int_{K^*} < x, y>^2   dx dy  \le \dfrac{n \omega^2}{ (n+2)^2}
  \end{equation}

The proof of this in [10] is self contained    using \\

\underline{Lemma5.2}([10] Lemma2.1) due to earlier authors; see reference in [9]

Let $f_1, f_2, f_3  :R^n \longrightarrow R_+$  be unconditional measurable  and satisfy the following inequality:-

$f_1 (x_1,x_2, \dots x_n) f_2(  y_1, y_2, \dots y_n)  \le 
f_3 (\sqrt{x_1 y_1 }, \sqrt{x_2 y_2} \dots \sqrt{x_n y _n }   ) ^2$  \textbf{...( c1)}  
then  

$\int f_1  (x) dx \int f_2 (y) dy  \le  ( \int f_3(z) dz )^2 $   \textbf{(c2)}

Conditions on equality "IFF"    are as follows:---

there is a continuous function $\tilde{f_3}$  non negative valued defined i $R^n$  such that all functions satisfy  these . We have used  the \underline{abbreviations} 
$ x= ( x_1, x_2, ......x_n) $  and likewise for  y and also c as vectors in $R^n$   and the unconventional, $  1/c=( 1/c_1,  1/c_2,.......1/c_n)$ assuming no zeros in denominators . Also  " products"   $xc $ and $ y 1/c$ where 
cx =   $(c_{k} x_{k}  :k=1, 2, ...n)  $  and    we define 
$y 1/c$ similarly

(i)$ \tilde{f_3} = f_3  $  ae and $\tilde{f_3}( x)  \tilde{f_3} (y)  \le \tilde{f_3 } ( \sqrt {x_1  y_1},  \sqrt{x_2 y_2} ............ \sqrt {x_n y_n}   )  $

(ii) there are constants  $ c_k$  for   $1\le k \le n $ and $ d>0$    such that

$  f_1 ( x) = d  \tilde{ f_3 }( c x) $ and $ f_2(y) =   1/d\tilde {f_3 }(  y 1/c) $  \\

\underline{The Proof of Th.5.1   ( outline)} 

Roughly, the idea is  : using coordinates expand the inner product in Left Hand side.Now etsimate individual term, there are "n" equal such. Let us give again only some details

We  may assume that the  usual vector basis $e_1, e_2,  \dots  e_n$ lies on boundaries of both  K  and   K* thanks to  invariance of Left Side in Ball's Inequality above  under  maps  from GL(n).

The object is  to define 3 functions to be able to apply Lemma5.3:--
First  let us  recall  $ e_1 ^ \perp $ the   hyperplane thru O and perpendicular to the vector $e_1$.
Now we  define these functions ( each from R to R)  using  volumes of parallel sections of  K   , or  K*  thus:--

\begin{equation}\label{E:conj}
f_ 1 (r) = | ( e_1 ^\perp + re_1) \cap K)|   and 
f_2(r) =   | ( e_1 ^\perp + re_1) \cap K^*)|   
\end{equation}

We note that both these  functions    are zero for $ r >1$  .Next suppose  that r and s are both in [0, 1]  with
$ f_1(r)  , f_2(s)  >0$

Now define 2 convex sets by
D(r)=$  e_1^\perp  \cap  ( K -  r e_1)  $ and  E(s)   using   s  and K*  in place of  r  an  K. Since these are unconditional, we have following:--

if	$u\epsilon D(r) $ ,  $v\epsilon E(s)$  then 
$\pm u $ is in D(r) and likewise  $ \pm v $  is in E(s) and so   $ r e_1 \pm u$  is in D(r)  and   $sv_1 \pm v $ is in E(s).	
The outcome is that the set   D(r)  in the space $e_1^\perp$  of dimension (n-1) and $D(r)  ^ *$ in the dual , and we have the inclusion

$ E(s)  \subset   ( 1-rs)   D(r)^* $   and  we obtain

$ |D(r)| | E( s)|  \le ( 1-rs) ^{n-1}  |D(r)| | D(r)^ *| 
\le  (1-rs)^{n-1}  \omega_{n-1}^ 2 $  ; this last due to the Blashcke Santalo Inequality in n-1 dimensions .

Observing from earlier  definitions that $  f_1(r) = |D(r) | $   and  likewise another   we get from the last 
$ f_1(r ) f_2(s)  \le   (1-rs)^{n-1} \omega_{n-1} ^2 $
and now we are in a position to define our 3rd   function:

\begin{equation}\label{E:conj}
f_3( t)  =t^2  ( 1-t^2) ^ {n-1/2} if  0\le t \le 1  and=0  if not.
\end{equation}

The  last inequality  enables us to apply Lemma 5.3   with n=1   due to  

$ f_1(r)r^2 f_2(s) s^2 \le  f_3( \sqrt{ rs} )  ^ 2  $ 
which   is \textbf{(c1)} 

Thus, we may use (c2) in this case    and obtain:-

( all integrals are over   $R+$  due to definition  of these functions)

$ \int f_1(r) r^2 dr  \int f_2(s) s^2 ds \le (\int f_3(t)  dt )^2 $

Recalling definitions of  $f_1  , f_2$  we   see the first integral  is $ 1/2 \int_K  x_1 ^2  dx$ and similarly with the second on the left  over K*    in  place of K.
The integral in Right  within parenthesis is  ( after a calculation with Gamma integral)  equal to 
$  \frac{\omega_n}{ 2( n+2)} $.

Thus all of this last translates into

$ 1/4\int_{K} \int _{ K^*}  x_i^2 y _i ^2   dx dy 
\le \frac {\omega_n ^2} { 4( n+2)^2} $, 

This estimate will be used below as one   of n equal terms

Further, the fact that  K is unconditional has as a result  

$\int_{K} x_i  x_j   dx =0 ; if i \ne j$

Now we  may   finish    proof of  Ball'Inequality in Statement:
we  have:--

\begin{equation}\label{E:conj}
\int_{K} \int_{K^*}  < x , y> ^2 dx dy 
= \Sigma_{i=1}^{ n}  \int_{K} \int_{K^*}  x_i ^2y_i^2  dx dy 
\end{equation}

  (  this upon expanding  scalar product , the product terms   vanish as noted above)  

Hence we  may use the  estimate for one   term  above  and  so the last is  

$\le  \frac{n\omega_n ^2}{  ( n + 2 )^2}  $  which is what we   wanted to prove.
Hence Ball's Inequality. QED \\

 For  the \underline{equality  statement}.  We  note that the hypotheses , for n=1 , in Lemma 5.3 are satisfied by the   eqs. above involving $f_1(r) r^2$ and likewise the others.
We may hence draw its conclusions:--

$f_1(r) r^2  =df_3(cr)$  and $f_2(r)r^2)= \frac{1}{d} f_3( r/d)  $ 

where we do not need " a.e" as the functions are continuous.  Also  we observe that c=1 , and   hence we  evaluate these functions:--

$f_1(r) =  d ( 1  -  r^2) ^{n-1/2}   \omega_{n-1} $  and

$f_2(r)  =\frac{1}{d} ( 1 - r^2) ^{ n-1/2} \omega_{n-1}$
Thereby we have

$f_1(r) f_2(r)  =  |D(r)| |E( r)| = (1-r^2)^{n-1} \omega_{n-1}^2$
and this yields equality   in earlier  inequality  with  both r and  s.
Therefore  we may apply te equality in Blaschke -Santalo Inequality (eq1)  for  n-1 .Thus   D(r) is an ellipsoid in dimension (n-1)   for  every $ r \epsilon [0, 1] $.
We\underline{  claim}  that in  fact, D(r) is an n-1  dimensional ball with radius   $ ( 1-r^2 ) ^ {1/2}$

Also  the unit vectors  $ e_k : k=1, 2, .....n $  are in the boundary of K ,and  as  D(0)   is an unconditional ellipsoid  . Hence  D(0)  is in fact the  n-1 dimensional unit ball in the space$ e_1^\perp $ . For the rest of the is  $i\ge2$  denote  $e_i^\perp  \cap (K-r e_i )  =  D_i(r)$
As above we may conclude that the set  $D_i(0)$  is an n-1  dimensional unit  ball in $e_i ^\perp$  for all such i  s.Hence  $D(r) $ is an unconditional ellipsoid and passing  thru the   n-1 points  $ a_2 e_2  ,   a_3e_3 ........a_n  e_n$    with each $x_k= \pm (1- r^2) ^{1/2}$  thereby proving our claim.

With the claim in hand we  may at last   conclude that K is  indeed the unit ball in $R^n$  . Since the Left side in Ball's Inequality is preserved under ( invertible) linear  maps of K we reached the desired end   in equality case QED \\
 
  \underline{{Examples }}   Let us illustrate  Ball's Inequality (27)   .In all these e the   product terms   integrate to zero as we saw in the proosf  . Hence we are left with square temrs of which  one is   $x_j^2   y_j^2$. Let  us recall $x=  (x_1, x_2........x_n)   $   and similarly for y

  \underline{Remark:}  The  equality when true, is  for each term   separately  in the examples we see below.   \underline{We dont know }if this is true in general\\

 Ex 4A    Let us consider  the Euclidean closed ubit disk in 2 dimensions, $K=K^0  = B_2$ As a  Calculus Exercise   we may verify the two sides in   Inequality (30) are equal f as  we indicate  briefly  As noted above    we  are left with square terms .  Consider one of them, $x_1  ^2  y_1^2$  and  $\int_{B} x_1 ^2  dx$  $(y_1^2)$.
 This integral in familiar terms  is  the double  integral   $\int_{B}x^2  dx dy$ .We may use polar coordinates to compute this  one  to be  $\pi /4$. The other integration  of $ y_1^2$ Over "B*  =B   also  gives same $\pi/4$ . Thus one of the square terms  we  have 
 
 $ \int_{B} \int_{B^*}  x_1  ^2  y_1 ^2  dx dy =   \pi ^2  /16 $
 Finally w have for Left Hand Side (LHS)    = $ 2 \pi^2 /16   =   \pi^2/8$.  The  RHS is  $ \frac{2} { 4^2}  \pi^2 $
 =$\pi^2  /8$ also , completing the example\\
 
 \underline{ Ex 4B}  This is of negative nature  with strict Inequality. Let us consider the  Unit  Square  $B_\infty $  or  B  ie   defined by  $ |x_1| ,|x_2|\le1   $  THe polar is the diamond  D ie unit ball in $l_1$ 
 
 The RHS   is  $\pi^2   /8$   . THe Left needs  calculation as in Ez2  ; there are 2 squae terms, each works  out to  $4/9$    and so LHS equals  $8/9$  $< $  RHS

\underline{ Remark}:  As   noted earlier in Introduction  that the calculations for the Integral or Left Hand side of Ball's   "quadratic" Inequality  Inequ (30)  are  harder than the ones in Exs of SEc 3.2. Let us only mention

\underline{Ex4B}

 .\\

  Following  [1], [2], [10]  we discuss the corresponding generalizations of  Blaschke-Santalo Inequality  eq(1) .  One example. of  it,  due to earlier authors  ( see  [1],[2] [9]and its references ) is the following  Th  5.1.  we recall  Eq(26)  the polar $f^o$  of a non negative valued ( measurable) function f on $R ^n $  . ( In  following all integrals are over $R^n $):- \\

\underline{ 5.3 Theorem} ( [1] [2] [10]Inequality 1.2)

Let f be an even, non   negative  Integrable function ;then

\begin{equation}\label{E:cong}
	\int f(x) dx  \int f^0  (x) dx  \le  ( 2\pi)^n
	\end{equation}
	with equality  iff f (x)= c exp  (- $ |Tx| ^2 $ )   for some positive definite matrix T and c $  >$0 .
	
	For a proof with some conditions on the function we refer to [1].	  
	
	 As  noted above in   [1] [2]  ( see  other references in [10])   this does reduce to eq(1)  Blaschke-Santalo   Inequality as we see below (in Remark 2) .\\
	 
	 \underline{ Ex 6 }Let us  note in contrast  that we   have strict Ineq  in (31)   with the characteristic function of  the unit Euclidean Ball  K in $R^2$ in   Ineq (20) . For  its polar, we need Fact1 above; using it ,  a Calculus 3 exercise leads  us  to:--  Left  Hand Side= $2\pi^2  <   4\pi ^2$ =Right Hand  Side \\

	 \underline{Remark 2:}  Let  us   reproduce  from [1] [3]the following:
	 
	 \underline{Ineq(31) implies  the upper bound  in Blaschke-Santalo Inequality}
	 
	 Given our convex body K, we recall from Introduction  the  gauge function $||x||_ K$  and let us denote it by  just $|| x || $  and likewise for the polar K* . The implication depends on the following formula  (  see reference in [1]):--
	 
	 $\int exp (- ||x||^2   /  2 ) dx =  c_n vol K  $
	 
	 where $c_n = \frac{(2\pi) ^ {n/2}}{vol B_n} $.
	 
	 with $B_n$ the Euclidean Ball.  We  apply  this  formula twice to the Left side  in  Ineq.  (20);  upon canceling the factor $( 2\pi) ^ n$  from both sides we have,
 
         $\frac{Vol K   vol K^*}{  (Vol B_n )^2}  \le 1 $
         and this gives the  Upper bound in  Blaschke- Santalo Inequality (1)\\

 Next  let us mention from [ 10 ]  a functional version of Ball's Inequality  and implies it .
      
   For this we  recall  the definition of "Unconditional  function"   above
 
 \underline{ 5.4 Theorem [ 10] Th  1.1}  

 Let f be  $\ge 0  $ integrable  on $R^n$  and  unconditional;  then we have

\begin{equation}\label{E:cong}
\int  \int (x,y) ^2   f(x)  f^0( y) dx dy   \le n ( 2\pi)^n
\end{equation}

There are iff conditions on equality above; these are same as the one for Theorem 5.1 except  that now T needs also to be diagonal 
 
 Using the same function   as in Th 5.1  we deduce Ball's Inequality   Th 5.2 from this one. 
  
  The proof  in [10] is also involved  and depends on idea   and results  on  log concave function and Ths 3.1, 3.2 there.

	\underline{MR Classifications:} Primary 52A20,  52A40, Secondary 42A05

\underline{	key words and phrases }; Blaschke- Santalo inequality, zonoid, unconditional bodies
	
	Retired from: SUNY /College;    Old Westbury NY 11568-0210
	rajan.anantharaman@gmail.com\\
	
	\textbf{Appendix}
	More recent reference:

	L.Grafakos,  Classical Fourier Analysis Grad Texts in Math ;   2014 Springer-Verlag,  , Heidelberg-New York\\

	\begin{tabular}{|l| r|r| r|r  |r |r |r|r |r |} \hline
		Name  &   $  $   &   &   &   &  &  &    \\  \hline
		Set   & A & $\subset$ &B &$ \longrightarrow$ & $ B  ^ o $ & $\subset$ & $ A  ^o $	\\	\hline
		function &f   & $\le$ &g & $\longrightarrow$  & $ g^o$ & $\le$ & $f^o$ \\  \hline
		both  & $A^{oo}$ & = & A &  & $ f^{oo}$ & =& f\\ \hline
	\end{tabular}

	\end{document}